\documentclass[english,11pt]{amsart}
\usepackage{a4wide}
\usepackage{amsmath,amssymb,graphicx}
\usepackage{mathrsfs}
\usepackage{mathtools}
\usepackage{float}
\usepackage{dsfont}
\usepackage[all,cmtip]{xy}

\usepackage{tikz}
\usetikzlibrary{intersections,calc,arrows.meta, matrix}
\usetikzlibrary{arrows}
\usepackage{tikz}
   \usetikzlibrary{positioning}
\usepackage{bbding,wasysym}

\graphicspath{{figures/}}
\usepackage{hyperref}
\hypersetup{
	colorlinks=true,
	linkcolor=blue,    
}

\usepackage{frcursive}
\usepackage[utf8]{inputenc} 
\usepackage{amscd}

\def\di{\displaystyle}

\newtheorem{theorem}{Theorem}
\newtheorem{definition}{Definition}
\newtheorem{lemma}{Lemma}

\newcommand{\N}{\mathbb{N}}

\newcommand{\R}{\mathbb{R}}

\newcommand{\T}{\mathbb{T}}

\begin{document}
\title[Wong-Zakai Variational integrators]{Dynamics of Stochastic Hamiltonian systems I - Wong-Zakai Variational integrators}
\author{Jacky CRESSON$^1$ and Rouba SAFI$^{1,2}$}

\begin{abstract}
We construct variational integrators to study the dynamics of stochastic Hamiltonian systems using their Wong-Zakai approximation. This approach can be used to interpret and justify previous work of Wang et al (L. Wang, J. Hong,  R. Scherer, F. Bai. Dynamics and variational Integrators of stochastic Hamiltonian systems. International Journal of Numerical Analysis and Modeling, 6(4), 2009) on stochastic variational integrators. Numerical examples are provided.
\end{abstract}

\maketitle

$^1$ Laboratoire de mathématiques et leurs applications, UMR CNRS 5142, Université de Pau et des Pays de l'Adour-E2S.

$^2$ LaMA, Laboratoire de math\'ematiques et applications (Tripoli, Liban), Universit\'e libanaise. 

%\tableofcontents

\section*{Introduction}

Stochastic Hamiltonian systems or {\bf Hamiltonian diffusion} are stochastic differential equations in the sense of Stratonovich \cite{oksendal} first introduced by J-M. Bismut in his seminal work called "Mécanique aléatoire" (Random Mechanics) \cite{bismut}. Hamiltonian diffusions consist in a drift part which coincides with a classical deterministic Hamiltonian system and a diffusion part which is a sum of Hamiltonian like quantities. Precisely, stochastic differential equations are called a {\bf stochastic Hamiltonian system} if it can be written in the form :
\begin{equation}
\label{stochastic-sys}
 \left\lbrace 
 \begin{array}{ll}
dp =& - \di\frac{\partial H}{\partial q}dt - \di\sum_{k=1}^{m} \di\frac{\partial H_k}{\partial q}\circ dW_k, \\
dq =&  \di\frac{\partial H}{\partial p}dt + \di\sum_{k=1}^{m} \di\frac{\partial H_k}{\partial p}\circ dW_k,
\end{array}
\right.
\end{equation}
where the $W_k$, $k=1,\dots ,m$, are independent Wiener process, $H$ and $\{ H_k \}_{k=1, \dots, m.}$ are functions from $\R^{2d}$ to $\R$. Such a definition extends easily to (symplectic) manifolds (see \cite{cami}).\\

The terminology of stochastic Hamiltonian systems is not clear since the function $H$ (or any $H_k$, $k=1,\dots ,m$) does not corresponds to the total energy of the system and is more delicate to interpret. However, specific properties of classical (deterministic) Hamiltonian systems are preserved:

\begin{itemize}
\item Hamiltonian systems possess a {\bf variational structure} meaning that their solutions corresponds to critical points of a functional (Hamilton's principle).   

\item An important geometric property of Hamiltonian system is that the associated flow is {\bf symplectic} \cite{ratiu,arnold}, meaning that it preserves the canonical two-form $\omega =dp\wedge dq$. This property has strong qualitative consequences. In particular, the flow is volume preserving (Liouville's theorem) \cite{hairer}. 
\end{itemize}

These two properties are preserved by the stochastic version of Hamiltonian systems (see \cite{bismut,cami}).\\

Like in the deterministic case, it is difficult to obtain explicitly the solutions of stochastic Hamiltonian systems, and numerical methods are needed. However, due to the specific symplectic structure of these systems, classical numerical methods usually behave poorly under long-time numerical simulations. Specific numerical methods called {\bf symplectic integrators} have been designed to construct numerical integrators that preserve the symplectic character of the flow. In other words, if the numerical method is defined by a mapping 
$(p_{n+1} ,q_{n+1} ) =\Phi (p_n , q_n )$ then we require that this mapping is symplectic, ensuring that we have
\begin{equation}
    dp_{n+1} \wedge dq_{n+1} = dp_n \wedge dq_n,\;\; \forall n \geq 1.
\end{equation}
In the deterministic case, previous work in this direction were made in particular by Vogelaere \cite{vogelaere}, Feng Kang et al. \cite{feng1},\cite{feng2},\cite{feng3}, etc. More references can be found in the classical textbook of E. Lubisch, Wanner and Hairer \cite{hairer} about geometrical numerical integrators. A way to construct symplectic integrators was developed by J.E. Marsden and co-workers by focusing on the discretization of the underlying functional. They are called {\bf Variational integrators}. A review can be found in the article of J-E. Marsden and M. West \cite{marsden}. \\ 

In the stochastic case, symplectic integrators were constructed in particular by G.N. Milstein in a series of papers (for example \cite{Milstein,Milstein2}) using a {\bf direct methods} meaning that he first derives conditions under which a general numerical scheme preserves symplecticity and then construct integrators. Using the Marsden's approach, two a priori different constructions of variational integrators were constructed by N. Bou-Rabee in  \cite{bou-rabee} and by L. Wang and co-workers in \cite{wang}): 

\begin{itemize}
    \item The construction of L. Wang et al. \cite{wang} is difficult to follow because it identifies two classes of objects of different nature. In particular, they identify the Stratonovich integral $\di\int_a^b f(t) \circ dW_t$ with $\di\int_a^b f(t) \dot{W} (t) dt$ although $W$ is nowhere differentiable so that all the computations and integrators developed in \cite{wang} need to be discussed more closely or amended. This is in particular the case for the variational formulation and proof proposed in \cite{wang} (see Section \ref{WZ}).

    \item In \cite{bou-rabee}, the method follows Marsden's approach to variational integrators by using classical approximation of the Stratonovich integral. However, they consider a special class of stochastic Hamiltonian systems where the configuration variable $q$ of every solutions is differentiable with respect to $t$, or equivalently that all the $H_k$ does not depend on $p$.
\end{itemize}

In this article, we compare the two constructions.\\

First, we give a meaning to the idea of L. Wang et al. \cite{wang} using the classical notion of {\bf Wong-Zakai diffusion approximation} \cite{zakai1,zakai2}. Formally, a smooth approximation $W_{\epsilon}$ of $W$ is constructed so that $\lim_{\epsilon \rightarrow 0} W_{\epsilon} =W$. As an example, one can use the classical averaging
$$W_{\epsilon} (t) =\di\frac{1}{\epsilon} \di\int_t^{t+\epsilon} W_s ds.$$
Using this approximation, a Stratonovich differential equations
$$dX=a(t,X)dt+\sigma (t,X) \circ dW_t ,$$
is obtained as the limit when $\epsilon$ goes to zero of the family of {\bf random differential equations}
$$\di\frac{dX_{\epsilon}}{dt} = a(t,X_{\epsilon} ) + \sigma (t,X) \dot{W_{\epsilon}},$$
the randomness coming from the term $\dot{W_{\epsilon}}$.\\

Applying this method to stochastic Hamiltonian system \eqref{stochastic-sys}, we obtain a one parameter family of non-autonomous (random) Hamiltonian systems defined by the Hamiltonian $H_{\epsilon}$ given by 
$$H_{W-Z,\epsilon} = H_0 +\di\sum_{k=1}^m H_k \dot{W_{\epsilon ,k}}$$
converging to the stochastic Hamiltonian system \eqref{stochastic-sys}. It must be noted that Wong-Zakai approximation for stochastic Hamiltonian systems has already been used at the beginning by J-M. Bismut \cite{bismut} in his study of the properties of Hamiltonian diffusions (see \cite{bismut}, Chapter 1 p.36 and Chapter 5, Section 2 p.224). Of course, one has to be careful because, as already noted by J-M. Bismut (\cite{bismut},p.26-27) not all the properties of classical Hamiltonian can pass to the limit. \\

These Hamiltonian systems can of course be obtained via the Hamilton principle by minimising the random functional 
$$\mathscr{L}_{H_{WZ,\epsilon}} (p,q)=\di\int_a^b (p\dot{q} -H_{WZ,\epsilon} (p,q,t)  ) dt$$
The idea is then to use well-known techniques for the construction of variational integrators of deterministic Hamiltonian systems to construct variational integrators for stochastic Hamiltonian systems.\\

We then have the following diagram:
$$
\xymatrix{
    X \ar[d]^{v.p}
    & 
    \ar[l]^{\epsilon \rightarrow 0} X_{\epsilon} \ar[d]^{v.p.}
    \\
    \mathscr{L}_{\mathbf{H}}  & 
    \ar[l]^{\epsilon \rightarrow 0} \mathscr{L}_{H_{\epsilon}} \\
    \ar[u]^{h\rightarrow 0}\mathscr{L}_{\mathbf{H},h} \ar[d]^{d.v.p.} & \ar[u]^{h\rightarrow 0} \mathscr{L}_{H_{\epsilon},h} \ar[d]^{d.v.p.}\\
    X_n  & \ar[l]^{ ? \epsilon \rightarrow 0} X_{\epsilon ,n}
}
$$
where {\it v.p.} and {\it d.v.p.} stand for the Hamilton principle and discrete Hamilton principle respectively. $X_n$ and $X_{\epsilon ,n}$ represent the variational integrators obtained from $\mathscr{L}_{\mathbf{H}}$ and $\mathscr{L}_{H_{\epsilon}}$ respectively.\\

It must be noted that the convergence of $X_{\epsilon,n}$ to $X_n$ is not trivial and is related to the existence of a limiting scheme consistent with the limiting equation. Such scheme are called notion {\it asymptotically preserving scheme} in \cite{bre1,bre2}. \\

In particular, the result depends drastically on the quadrature used in the Wong-Zakai approximation as explained in \cite{londono}. As an example, let us consider the linear Stratonovich stochastic differential equation
\begin{equation}
    dX_t =a X_t dt +b X_t \circ dW_t .
\end{equation}
For $X_0 \in \R$, the solution is given by $X_t = X_0 \exp (at+b W_t )$. The Wong-Zakai approximation leads to the one parameter family of random differential equations
\begin{equation}
    \di\frac{dx}{dt} = ax +bx \dot{W}_{t,\epsilon} .
\end{equation}
As we have an ordinary differential equation, one can use all the classical tools to discretize it. In figure \ref{fig1}, we present the simulations $X_1$ obtained using the Euler scheme, $X_{21}$ mixing a mid-point formula for the deterministic part and a Euler scheme for the random part, $X_{12}$ mixing a Euler scheme for the deterministic part and a mid-point formula for the random part, $X_2$ where a mid-point formula is used for both part and finally the exact solution denoted $X_E$. We have taken $a=1.5$ and $b=1$, $N=100$ and the increment $h=0.01$.

\begin{figure}[!ht]
    \centering
    \includegraphics[width=0.8\textwidth]{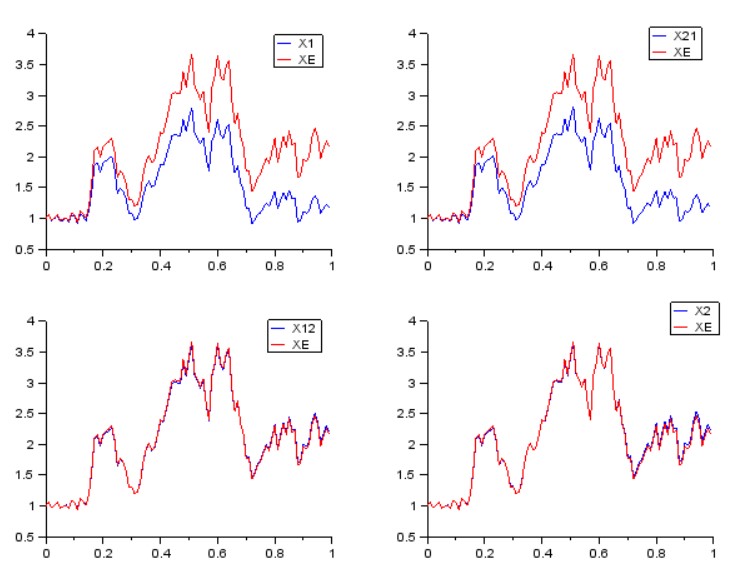}
    \caption{Different discretizations of the Wong-Zakai approximation}
    \label{fig0}
\end{figure}
We observe in figure \ref{fig0} a good agreement between the exact solution and simulations as long as a quadrature of order $2$ is used for the random part. This phenomenon is ultimately related to the definition of the Stratonovich integral which use a mid-point quadrature formula in its definition via Riemann sums which is of order $2$.\\

In this article, we construct the Wong-Zakai variational integrator and the stochastic variational integrators using the {\bf discrete embedding} theory (see \cite{cr1,cr2,cr3}), i.e. by defining a discrete differential and integral calculus and the corresponding calculus of variations. In order to cover Stratonovich stochastic differential equations, we use the mid-point embedding developed in \cite{rouba1,rouba2} in the deterministic case. This approach allows us to obtain a new derivation of the Bou-Rabee et al. stochastic variational integrators \cite{bou-rabee} and a rigorous derivation of Wong-Zakai variational integrators as a justification for computations made by Wang et al. in \cite{wang}.\\

Our article is organized as follows: \\

%Section \ref{stochamilton} remind the classical definition of stochastic Hamiltonian systems and there properties. 
Section \ref{midpointemb} remind the mid-point embedding formalism developed in \cite{rouba2}. We use this formalism in Section  to define a discrete stochastic calculus of variations. This allows us to define in Section \ref{WZ-VI} a notion of discrete stochastic Hamiltonian system which coincide in the deterministic case to the Marsden-Wendlandt definition of discrete Hamiltonian system \cite{wendlandt}. section \ref{example}
\section{Wong-Zakai approximation Hamiltonian}
\label{WZ}

Let $S_{\mathbf{H}}$ be a stochastic Hamiltonian system defined by the family $\mathbf{H}=\{ H_0, \dots ,H_k \}$. We denote by $W_{i,\epsilon}$, $i=1,\dots ,k$ a smooth approximation of the $k$ independent Wiener processes $W_i$, $i=1,\dots ,k$.\\

Let us define the {\bf Wong-Zakai approximation Hamiltonian} (or simply Wong-Zakai Hamiltonian) of the stochastic Hamiltonian system $S_{\mathbf{H}}$ by: 
\begin{equation}
H_{W-Z,\epsilon} (p,q,t)=H_0 (p,q)+\di\sum_{i=1}^k H_i (p,q) W_{\epsilon ,i} (t,\omega ) ,
\end{equation}
for all $\epsilon >0$ and $\omega \in \Omega$.\\

Two remarks about the Wong-Zakai Hamiltonian:\\

\begin{itemize}
\item The Wong-Zakai Hamiltonian define a {\bf random Hamiltonian systems} in the sense of {\it random differential equations} \cite{han}, the randomness coming from the approximation of the stochastic part. 

\item Even if we begin with a family autonomous Hamiltonian system $(H, H_k, k=1,\dots ,m)$ the resulting Wong-Zakai Hamiltonian is non-autonomous. As a consequence, these Hamiltonian represent in general nonconservative dynamics.  
\end{itemize}
\vskip 2mm
The previous remarks allows us to interpret which kind of dynamics stochastic Hamiltonian systems model: Taking $H$ as an initial Hamiltonian dynamics which is conservative, one consider non-autonomous Hamiltonian random perturbations of $H$ leading to non-conservative random Hamiltonian dynamics. The stochastic Hamiltonian dynamics corresponds to the asymptotic dynamics generated by these non-conservative random Hamiltonian dynamics. The randomness implies that we can have dissipation or increasing of energy. This point of view can be compared with the one presented in Wang et al. in \cite{wang} which takes a reverse presentation: they begin with modeling of nonconservative force and say that stochastic Hamiltonian systems can be seen as "Hamiltonian systems ... disturbed by certain nonconservative force" called "random force" in (\cite{wang}, p.589).\\

A natural question with respect to stochastic Hamiltonian systems and their Wong-Zakai approximations is related to the behaviour of energy. Indeed, when no stochastic perturbation is present, we have
\begin{equation}
    \di\frac{d}{dt} (H(p_t ,q_t ,t))= \di\frac{\partial H}{\partial t} (p_t ,q_t ,t) ,
\end{equation}
which implies for autonomous Hamiltonian systems $H(p,q)$ the conservation of energy represented by $H$.\\

In the stochastic case, even the formulation of such a conservation property is difficult because we have not a single function representing energy for such systems. contrarily to the case of the Wong-Zakai approximations where $H_{WZ,\epsilon}$ stands for the energy of the system. Two questions can be studied:\\

\begin{itemize}
    \item The only object having some intrinsic meaning of energy in the stochastic case is the unperturbed energy $H$. As a consequence, we can look for the preservation of $H$ under the dynamics of the stochastic Hamiltonian system.

    \item Does the behaviour of the energy $H_{W-Z,\epsilon}$ say something about the dynamics of the stochastic Hamiltonian system ?

    \item Do we have connection between the first integrals of the stochastic Hamiltonian system and its Wong-Zakai approximations ?
\end{itemize}
\vskip 2mm
The answer is not so simple. \\

As an example, if one consider the Kubo stochastic Hamiltonian system
\begin{align}
    \label{kubo1}
    dp &= -aqdt - \sigma q \circ dW_t, & p(0)=p_0,\\
    \label{kubo2}
    dq &= apdt + \sigma p \circ dW_t, & q(0)=q_0,
\end{align}
with
\begin{align}
    H(p, q) = \di \frac{a}{2}(p^2 + q^2),&&  H_1(p, q) = \di \frac{\sigma}{2}(p^2 + q^2),
\end{align}
we can prove the following simple result:

\begin{lemma}
    \label{kuboener}
    For the stochastic Kubo oscillator we have:
    \begin{itemize}
    \item The unperturbed energy $H$ is preserved under the stochastic perturbation.
    
    \item All the Wong-Zakai Hamiltonians $H_{W_Z,\epsilon}$ are nonconservative.

    \item The unperturbed energy $H$ is preserved under all the dynamics of the Wong-Zakai approximations.
    \end{itemize}
\end{lemma}

Indeed, for all $\epsilon >0$, we have as usual
\begin{equation}
\di\frac{d}{dt} \left ( H_{WZ,\epsilon} (p_t ,q_t ) \right ) = \di\frac{\partial H_{WZ,\epsilon}}{\partial t} (p_t ,q_t ) = \di\sum_{i=1}^{m} H_k (p_t ,q_t) W"_{\epsilon,k} (\omega ) , 
\end{equation}
if each $W_{\epsilon,k}$ is a sufficiently smooth approximation of the Wiener process $W_k$. As $W"_{\epsilon ,k}$ will take arbitrary positive or negative values, we will obtain a nonconservative system.\\

However, looking for the behaviour of $H(p_t ,q_t )$ over the solution of the Wong-Zakai Hamiltonian, we obtain
\begin{equation}
\di\frac{d}{dt} \left ( H (p_t ,q_t ) \right ) = \di\sum_{k=1}^m \left ( \di\frac{\partial H}{\partial p} \di\frac{\partial H_k}{\partial q} -\di\frac{\partial H}{\partial q} \di\frac{\partial H_k}{\partial p} \right ) 
\dot{W}_{\epsilon ,k} =\di\sum_{k=1}^m \{ H,H_k \}  \dot{W}_{\epsilon ,k}
\end{equation}
where $\{ \cdot ,\cdot \}$ is the classical Poisson bracket defined for two functions by $\left \{ f,g\right \} = \di\frac{\partial f}{\partial p}\di\frac{\partial g}{\partial q} -\frac{\partial f}{\partial q}\di\frac{\partial g}{\partial p}$.\\

In the Kubo case, we have $H_1 =(\sigma /a) H$ so that $\{ H,H_1 \}=0$ and $H$ is a first integral of the Wong-Zakai dynamics.\\

This property goes to the limit. Indeed, using the Stratonovich differential calculus, one have 
\begin{equation}
\label{energy}
d \left(H (p_t ,q_t ) \right)  = \di\sum_{k=1}^m \left ( \di\frac{\partial H}{\partial p} \di\frac{\partial H_k}{\partial q} -\di\frac{\partial H}{\partial q} \di\frac{\partial H_k}{\partial p} \right ) dW_k  =\di\sum_{k=1}^m \{ H,H_k \}  dW_k 
\end{equation}
so that $dH=0$ for the stochastic Kubo oscillator meaning that the stochastic Hamiltonian preserve the unperturbed energy $H$ of the system.\\

This simple example shows that the behaviour of the Wong-Zakai energy $H_{W-Z,\epsilon}$ does not give significant information about the dynamics of the stochastic Hamiltonian systems in the contrary to the unperturbed energy $H$.\\

Moreover, the use of stochastic Hamiltonian systems as model for "nonconservative" forces acting on the unperturbed Hamiltonian $H$ seems not so clear. Indeed, in the Kubo case, as proved, the unperturbed energy is preserved. However, the Kubo example is highly non generic. Indeed, using \eqref{energy}, it is easy to see that in order to ensure the preservation of the energy $H$ for the stochastic Hamiltonian system, one needs to have $\{ H,H_k \}=0$ for all $k=1,\dots ,m$  (see \cite{bismut}, Theorem 4.2 p.230). These commutation properties are in general not satisfied so that the unperturbed energy is usually not preserved. As a consequence, stochastic Hamiltonian systems are generically associated to "nonconservative" behaviours. \\

The computations made previously shows that:

\begin{lemma}
    If $I$ is a first integral of the stochastic Hamiltonian system then $I$ is a first integral of each Wong-Zakai Hamiltonian. 
\end{lemma}

As a consequence, Wong-Zakai Hamiltonians preserve many qualitative properties of the stochastic Hamiltonian system and are then good candidates for the construction of numerical integrator.\\ 

We denote by $S_{H_{\epsilon}}$ the associated Hamiltonian system defined by:

\begin{equation}
\label{WZham}
    \left .
    \begin{array}{lll}
    \di\frac{dp}{dt} & = & -\di\frac{\partial H}{\partial q} -\di\sum_{k=1}^{m} \frac{\partial H_k}{\partial q} \dot{W}_{\epsilon,k} ,\\
    \di\frac{dq}{dt} & = & \di\frac{\partial H}{\partial p} +\di\sum_{k=1}^{m}  \frac{\partial H_k}{\partial p} \dot{W}_{\epsilon,k} ,
    \end{array}
    \right .
\end{equation}

By classical result on Hamiltonian systems, we have:

\begin{theorem}
\label{theo1}
    Solutions of $S_{H_{W-Z,\epsilon}}$ correspond to the critical points of the random functional
\begin{equation}
\label{WZfunc}
\mathscr{L}_{H_{W-Z,\epsilon}} (p,q)=\di\int_a^b (p\dot{q}-H_{W-Z,\epsilon} (p,q,t)) dt .
\end{equation}
\end{theorem}

Theorem \ref{theo1} is a simple consequence of the classical Hamilton's principle for each Wong-Zakai Hamiltonians. It must be noted that the corresponding result in Wang et al. \cite{wang} is what they call the stochastic Hamilton's principle (see \cite{wang}, Theorem 2.3 p. 591) as long as what they call stochastic Hamiltonian system written as (see \cite{wang}, equations (18)-(19)):

\begin{equation}
\label{WZhamwang}
    \left .
    \begin{array}{lll}
    \di\frac{dp}{dt} & = & -\di\frac{\partial H}{\partial q} -\di\sum_{k=1}^{m} \frac {\partial H_k}{\partial q} \circ \dot{W}_{k} ,\\
    \di\frac{dq}{dt} & = & \di\frac{\partial H}{\partial p} +\di\sum_{k=1}^{m} \frac {\partial H_k}{\partial p} \circ \dot{W}_{k} ,
    \end{array}
    \right .
\end{equation} 
where $W_k$, $(k=1, \dots, m)$ is assumed to be independent Wiener processes. However, the previous system has no meaning from the mathematical point of view unless we interpret it as \eqref{WZham}. Indeed, if \eqref{WZhamwang} is to be interpreted as a Stratonovich saying $\circ \dot{W}_{t}$ stands for the notation $\circ dW_k$ then the left hand side must be $dp$ or $dq$ and a $dt$ must appear after the deterministic part. Moreover, as already pointed out by the authors (see \cite{wang}, after equation (17)) a Wiener process is nowhere differentiable so that the notation $\dot{W}$ is subject to caution. The main point is that the Lagrangian functional used by Wang et al. \cite{wang} whose critical points correspond to solutions of \eqref{WZhamwang} is given by 

\begin{equation}
\label{funcwang}
    \mathscr{L} (p,q)=\di\int_a^b (p\dot{q} -H(p,q)) dt -\di\sum_{k=1}^m \di\int_a^b H_k (p,q) \circ dW_{k,t} .
\end{equation}

However, this functional has no meaning if \eqref{WZhamwang} is interpreted as a Stratonovich equation unless $q$ is differentiable with respect to $t$ meaning that all the $H_k$ depend only on $q$. It must be noted that this assumption is precisely the one made by N. Bou-Rabee et al. in \cite{bou-rabee}. However, such assumption is not made in \cite{wang} so that \eqref{funcwang} must be interpreted as \eqref{WZfunc}. This point of view on \cite{wang} is reinforced by the fact that only the classical Dubois-Raymond theorem is used (see \cite{wang},Lemma 2.2 p.590) in the proof of (\cite{wang}, Theorem 2.3 p.591) using the fact that we have a term $\dot{W}_t dt$ in the functional and not a $\circ dW_t$.

\section{Discrete mid-point stochastic differential and integral calculus}
\label{midpointemb}

In this section, we recall the discrete mid-point differential and integral calculus developed in \cite{rouba2}. Using the definition of the Stratonovich integral, we then define the discrete stochastic integral in the sens of Stratonovich and we discuss some properties like a stochastic version of the Dubois-Raymond lemma.

\subsection{Mid-point differential and integral calculus}

We refer to \cite{bohner} for an introduction to {\bf time scale calculus} in general. We restrict our attention to uniform discrete time scale in the following. Time scale calculus is by nature an order one discrete integral or differential calculus. In order to keep the differential/integral presentation of time scale but to cover approximations of {\bf higher order} (here order two), we need to introduce several time scales: \\

Let $I=[a,b] \subset \R $,\, $N\in \N^*$ and $h=(b-a)/N$. We define the following time scales:
\begin{enumerate}
    \item[-]$\T=\{t_i=a+ih,\ i=0,1,...,N\}$.
    \item[-]$\T^+=\T \setminus \{b\},\;\;\T^-=\T \setminus \{a\}, \;\; \T^{\pm}= \T^+ \cap \T^- $. 
    \item[-]$\T_{\frac{1}{2}}=\{t_{i+\frac{1}{2}}=\frac{1}{2}(t_{i+1}+t_i),\;i=0,\dots,N-1\}$.
    \item[-]$\T_{\circ} =\T\cup\T_{\frac{1}{2}}$.
\end{enumerate}

For a given discrete time scale $\Tilde{\T}$ on $[a,b]$ with increment $\Tilde{h}$, we denote by $\sigma_{\Tilde{\T}}$ and $\rho_{\Tilde{\T}}$ the maps defined by 
\begin{align*}
\sigma_{\Tilde{\T}} \colon \Tilde{\T}^+  &\to \Tilde{\T}^- \\
t &\mapsto \sigma_{\Tilde{\T}}(t)=t+\Tilde{h},
\end{align*}
and 
\begin{align*}
\rho_{\Tilde{\T}} \colon \Tilde{\T}^-  &\to \Tilde{\T}^+ \\
t &\mapsto \rho_{\Tilde{\T}}(t)=t-\Tilde{h}.
\end{align*}

In the following, we use the simplified notations:
$$\sigma=\sigma_{\T}\;(resp.\,\rho=\rho_{\T}),\ \sigma_{\frac{1}{2}}=\sigma_{\T_\frac{1}{2}}\;(resp.\,\\ \rho_{\frac{1}{2}}=\rho_{\T_\frac{1}{2}}),\  \sigma_\circ=\sigma_{\T_\circ}\;(resp.\,\rho_\circ=\rho_{\T_\circ}).
$$

For all $f\in C(\T,\R^d)$, we define $f_\circ \in  C(\T_\circ,\R^d)$ as an extension of $f$ on $\T_\circ$  by
\begin{equation}
f_\circ(t)=\;
\left\lbrace
\begin{array}{ll}
f(t_i),& t=t_i.\\\\
\dfrac{f(t_i)+f(t_{i+1})}{2},& t=t_{i+\frac{1}{2}}.
\end{array}
\right.
\end{equation}

In the following, we use different operators. 

\begin{definition}
    Let $f\in C(\T_{\frac{1}{2}} ,\R^d )$. We denote by $[f]_{\frac{1}{2},-}$ the function defined on $\T_{\frac{1}{2}}^-$ by
\begin{equation}
    [f]_{\frac{1}{2},-} (t)=\di\frac{1}{2} 
    \left ( 
    f(t)+f (\rho_{\frac{1}{2}} (t))
    \right ) ,
    \ \ \mbox{\rm for all}\ \ t\in \T_\frac{1}{2}^- .
\end{equation}
Equivalently, for $f\in C(\T_{\circ} ,\R^d )$, we denote by $[f]_\circ$ the function defined on $C (\T^{\pm} ,\R^d )$ by 
\begin{equation}
    [f]_{\circ } (t) =\di\frac{1}{2} \left (
    f(\sigma_{\circ} (t) )+f( \rho_{\circ} (t) )
    \right ) ,
    \ \ \mbox{\rm for all}\ \ t\in \T^{\pm} .
\end{equation}
\end{definition}

For an arbitrary time scale $\tilde{\T}$ and a function $f$ of $C([a,b],\R^d )$ we denote by $\pi_{\tilde{\T}}$ the map from $C([a,b],\R^d )$ into $C(\tilde{\T},\R^d )$ obtained by taking the restriction of $f$ over $\tilde{\T}$.

\subsubsection{Discrete derivative and anti-derivative}

Let $\tilde{\T}$ be an arbitrary time scale. We denote by $\Delta_{\tilde{\T},+}$ and $\Delta_{\tilde{\T},-}$ the operators defined $for\ all\ f\in (\tilde{\T},\R^d)$ by 
\begin{equation}
\Delta_{\tilde{\T},+}f(t)=\di\frac{f^{\sigma_{\tilde{\T}}}(t)-f(t)}{ \sigma_{\tilde{\T}}(t)-t} ,  
\end{equation}
and
\begin{equation}
\Delta_{\tilde{\T},-}f(t)=\di\frac{f(t)-f^{\rho_{\tilde{\T}}}(t)}{t-\rho_{\tilde{\T}}(t)},   
\end{equation}
with 
\begin{equation}
f^{\sigma_{\tilde{\T}}} := f\circ \sigma_{\tilde{\T}}\ \mbox{\rm and}\ \  f^{\rho_{\tilde{\T}}} :=f\circ \rho_{\tilde{\T}} .
\end{equation}

According to the used time scale, we simplify our notations as follows:

\begin{equation}
\Delta_+ =\Delta_{\T,+} \ (resp.\; \Delta_-=\Delta_{\T,-} ),
\Delta_{\circ,+}=\Delta_{\T_\circ,+} \ (resp.\; \Delta_{\circ,-}=\Delta_{\T_\circ,-} ),
\Delta_{\frac{1}{2},+}=\Delta_{\T_\frac{1}{2},+} \ (resp.\; \Delta_{\frac{1}{2},-}=\Delta_{\T_{\frac{1}{2}},-} ).
\end{equation}

The associated discrete anti-derivative is then defined by (see \cite{rouba2}):

Let $\lambda \in [0,1[$. We denote by $t_{i,\lambda} =(1-\lambda )t_i +\lambda t_{i+1}$, $i=0,\dots ,N-1$. We denote by $\T_{\lambda}$ the set of $t_{i,\lambda}$, $i=0,\dots ,N-1$ and $\T_{\circ ,\lambda} =\T \cup \T_{\lambda}$. The $\lambda$-antiderivative over $\T$ is defined for all function $f\in C(\T_{\circ ,\lambda}, \R^d )$ by 
\begin{equation}
\int_{t_i}^{t_{i+1}} f(t) \Delta_{\lambda , \T} t= f (t_{\lambda}) (t_{i+1} -t_i ).
\end{equation}

It must be noted that despite the fact that we need the information about $\T_{\circ ,\lambda}$, the discrete antiderative is only defined on $\T$, meaning that we consider only integrals whose bounds of integration belong to $\T$.\\

We denote by $\iota_{0,\lambda ,\T,+}$ the mapping from $C(\T_{\lambda} ,\R^d )$ into $P^{0,+}_{\lambda,\T} ([a,b[,\R^d )$ defined by 
\begin{equation}
\iota_{0,\lambda ,+} [f] (t)=
f(t_{i,\lambda} )\  \mbox{\rm for all} \  t\in [t_i ,t_{i+1} [,\ i=0,\dots ,N-1.
\end{equation}

We simplify our notations according to the used time scale. Precisely, we denote by 
\begin{equation}    
\displaystyle\int_{a}^{b} f(t) \Delta t=\displaystyle\int_{a}^{b} f(t) \Delta_{\T} t,\ 
    \displaystyle\int_{a}^{b} f(t) \Delta_{\frac{1}{2}} t=\displaystyle\int_{a}^{b} f(t) \Delta_{\frac{1}{2},\T} t.
\end{equation}

\subsection{Properties of discrete derivatives and integrals}

\begin{lemma}
\label{relationprop}
For all $f\in C(\T,\R^d)$, we have 
\begin{equation}
\Delta_{\circ,+}[f_{\circ}](t_{i+\frac{1}{2}})=\Delta_+[f](t_i), \;\;\;\;\;for\, all\, i=0,...,N-1.
\end{equation}
\end{lemma}

\begin{lemma}[Discrete integration by part formula]
\label{itegration by part}
let $f\in C(\T_\circ,\R^d)$ and $v\in C(\T,\R^d)$, we have
\begin{equation}
\displaystyle\int_{a}^{b} f(t)\Delta_{\circ,+}[v_{\circ} ](t)  \Delta_{\frac{1}{2}} t=
-\displaystyle\int_{a}^{b} \Delta_{\frac{1}{2},-}[f]( \sigma_{\circ} (t) ) v(t)
\, \Delta t +f(t_{(N-1)+\frac{1}{2}})v(t_N)-f(t_{\frac{1}{2}})v(t_0)
.
\end{equation}
\end{lemma}

\begin{lemma}\label{dubois1}
For all $f\in C(\T_\frac{1}{2},\R^d)$, $v\in C_0 (\T,\R^d)$ we have 
\begin{equation}
\label{transfer}
\int_{a}^{b} f(t)v_\circ(t) \Delta_{\frac{1}{2}} t = 
    \int_{a}^{b} \left [ f \right ]_{\circ} (t) v(t)\Delta t + \frac{h}{2} \left( f(t_{N-\frac{1}{2}}) v(t_N)+ f(t_\frac{1}{2}) v(t_0) \right ) ,
\end{equation}
or equivalently 
\begin{equation}
\int_{a}^{b} f(t)\, v_{\circ} (t) \Delta_{\frac{1}{2}} t =  \int_{a}^{b} 
[ f ]_{\frac{1}{2},-} (\sigma_{\circ} (t)) v (t) \Delta t+ \frac{h}{2} \left( f(t_{N-\frac{1}{2}}) v(t_N)+ f(t_\frac{1}{2}) v(t_0) \right ).
\end{equation}
\end{lemma}

For more details see \cite{rouba2}.
\\

In order to conclude the Euler Lagrange equation in the continuous calculus, we apply the Dubois-Reymond lemma. One need to derive a discrete analogue of Dubois-Reymond lemma in the dicrete frame.

\begin{lemma}(discrete Dubois-Raymond lemma)
\label{stoc-duboi}
Let $f \in C(\T,\R^d)$ such that 
\begin{equation}
    \int_{a}^{b} f(t) v(t) \Delta t =0,\;\;\; for \ all \ v\in C(\T,\R^d),
\end{equation}
then we have  
$f(t)=0, \;\;\; for\ all \ t\in \T^{+,-}$.
\end{lemma}

\section{Wong-Zakai variational integrators}
\label{WZ-VI}
In this section, we derive variational integrators obtained using Wong-Zakai Hamiltonians and their variational formulations. As already reminded in the introduction, one can not use arbitrary quadrature formula for the discretization of the random part as it impacts the convergence of the resulting scheme to the solutions of the underlying stochastic Hamiltonian system. In order to illustrate this phenomenon, we compare two variational integrators obtained from the Wong-Zakai Hamiltonians using the strategy of discrete embedding exposed in previous work \cite{rouba1,rouba2}: an order one scheme and an order two which uses the mid-point embedding developed in \cite{rouba2}. We show that, as expected, the order one variational integrator does not reproduce the correct dynamics of the stochastic Hamiltonian system contrary to the mid-point scheme. This Section must be compared with the results exposed in (\cite{wang}, Section 2 p. 592-595). 

\subsection{Construction of discrete Lagrangian functional: principles}

The construction of the discrete functional is done using the discrete embedding formalism, i.e. thinking all the written objects using integrals and derivatives as certain integro-differential operators and constructing their discrete analogue using a given discrete differential and integral calculus.\\ 

The Hamiltonian functional for example can be understood as follows:\\

Let us denoted by $O (\di\int_a^b \cdot dt , \di\frac{d}{dt} ,H)$ the integro-differential operator acting on couple $(p,q)$ of $C^0\times C^1 $ functions by
\begin{equation}
O (\di\int_a^b \cdot dt ,\di\frac{d}{dt} ,H)(p,q):= \di\int_a^b (p\dot{q} -H(p,q,t) ) dt.
\end{equation}
By definition, the classical Hamiltonian functional $\mathscr{L}_H (p,q)$ is then defined by 
\begin{equation}
    \mathscr{L}_H (p,q)=O (\di\int_a^b \cdot dt ,\di\frac{d}{dt} ,H)(p,q) ,
\end{equation}
which indicates the specific role of the classical differential and integral calculus.\\

Let $\T$ be a given discrete time-scale with increment $h$ and let us denoted formally by $\Delta$ and $\di\int_a^b \cdot \Delta t$ the discrete differential and integral calculus defined over the functional space $C(\T ,\R^d )$ which is fixed. \\

The discrete functional $\mathscr{L}_{\T}$ is then defined for couple of functions $(p,q)\in C (\T^+ ,\R^d )\times C(\T ,\R^d )$ by 
\begin{equation}
\mathscr{L}_{\T} (p,q) :=O (\di\int_a^b \cdot \Delta t ,\di\Delta ,H)(p ,q ) . 
\end{equation}
The algebraic form of the integro-differential operator $O$ is exactly the same in the continuous and discrete case. The only difference is that the underlying integral and differential calculus is adapted to the functional context here to consider discrete functions. \\

Different choices in the discrete differential and integral calculus will lead to different discrete functional and naturally to different numerical scheme. It must be noted that some discrete operators need an extension of the functional space $C(\T ,\R^d )$ in order for example to obtain better quadrature formulas. An example of such phenomenon is given by the mid-point differential and integral calculus defined in the previous section. \\

This way to construct discrete analogue of functional can be compared with the one used by Wang et al. in (\cite{wang}, Section 2, p. 593). They fix first a set of time $t_i$ on a given interval $[a,b]$. This corresponds to the choice of a given discrete time-scale $\T$. The discrete functional, denoted by $\overline{\mathscr{S}}_h$, $h=t_{i+1}-t_i$, is defined on a finite family of couples $\{ (p_n ,q_n \}_0^N$ by choosing different quadrature formula for $\dot{q}$ and the integral of a function over $[t_i,t_{i+1}]$. The family $\{ (p_n ,q_n \}_0^N$ is the analogue of looking for two functions in $C(\T ,\R^d )$. The choice of the quadrature formula for $\dot{q}$ corresponds to the definition of a specific discrete differential calculus. A quadrature formula for a given integral over $[t_i ,t_{i+1}]$ corresponds to a specific discrete integral calculus. The basic ingredients can then be compared between the two approach. But, the discrete embedding point of view tells exactly how the integro-differential structure of the functional is changed by passing from a continuous setting to a discrete one. This property is completely lost in the usual construction proposed by Wang et al. in (\cite{wang}, (46) p.593).

\subsection{Examples of discrete functional}

Using the discrete embedding associated to the choice $(\Delta_+ ,\di\int_a^b \cdot \, \Delta t )$ we obtain the following order one discrete random functional:
\begin{equation}
    \mathscr{L}_{H_{W-Z,\epsilon},h,1} (p,q) :=O (\di\int_a^b \cdot \Delta t ,\di\Delta_{+} ,H)(p ,q ) , 
\end{equation}
or more explicitly 
\begin{equation}
    \mathscr{L}_{H_{W-Z,\epsilon},h} (p,q) := h \di\sum_{i=0}^{n-1} \left [ 
    p_i \left ( \di\frac{q_{i+1}-q_i}{h} \right ) - H_{W-Z,\epsilon} (p_i ,q_i ,t_i ) \right ] . 
\end{equation}

In the Kubo case, we obtain:
\begin{equation}
    \mathscr{L}_{H_{W-Z,\epsilon},h} (p,q) := h \di\sum_{i=0}^{n-1} \left [
p_i \left ( \di\frac{q_{i+1}-q_i}{h} \right ) -\di\frac{a}{2} \left ( p_i^2 +q_i^2 \right )
-\di\frac{\sigma }{2} \left ( p_i^2 +q_i^2 \right ) \Delta W_i \right ] ,
\end{equation}
where $\Delta W_i := W(t_{i+1})-W(t_i)$.\\

The mid-point embedding of the Hamiltonian functional used the mid-point differential and integral calculus defined in Section \ref{midpointemb}. We then obtain the discrete random functional $\mathscr{L}_{H_{W-Z,\epsilon},h}$ defined by 
\begin{equation}
\label{wong-zakai-discrete}
    \mathscr{L}_{H_{W-Z,\epsilon},h,2} (p,q) :=O (\di\int_a^b \cdot \Delta_{1/2} t ,\di\Delta_{\circ ,+} ,H)(p_{\circ},q_{\circ}) , 
\end{equation}
for all $(p,q)$ in $C (\T^+ ,\R^d )\times C(\T ,\R^d )$. \\

This case illustrates the fact that in order to obtain quadrature formula of higher order, one must enrich the discrete functional space $C(\T ,\R^d )$. Indeed if one want to define the mid-point quadrature formula from a given discrete function in $f\in C(\T ,\R^d )$, we have to consider an extension of $f$ to a function in $C(\T_{\circ} ,\R^d )$.\\

The mid-point discrete random functional is explicitly given by

\begin{equation}
    \mathscr{L}_{H_{W-Z,\epsilon},h} (p,q): = h \di\sum_{i=1}^{n-1} 
\left [ 
    \left ( \di\frac{p_i +p_{i+1}}{2} \right ) \left ( \di\frac{q_{i+1}-q_i}{h} \right ) - H_{W-Z,\epsilon} \left (  \di\frac{p_i +p_{i+1}}{2} , \di\frac{q_i +q_{i+1}}{2}, \di\frac{t_i +t_{i+1}}{2} \right ) 
\right ] . 
\end{equation}

In the Kubo case, we obtain:

\begin{equation}
\begin{array}{ll}
    \mathscr{L}_{H_{W-Z,\epsilon},h} (p,q) &= h \di\sum_{i=1}^{n-1} 
   \left [  \left ( \di\frac{p_i +p_{i+1}}{2} \right ) \left ( \di\frac{q_{i+1}-q_i}{h} \right ) - \di\frac{a}{2} \left ( \left (  \di\frac{p_i +p_{i+1}}{2} \right )^2 +\left (  \di\frac{q_i +q_{i+1}}{2} \right )^2 \right ) \right .\\
    &\left . -\di\frac{\sigma}{2} \left ( \left (  \di\frac{p_i +p_{i+1}}{2} \right )^2 +\left (  \di\frac{q_i +q_{i+1}}{2} \right )^2 \right ) \Delta W_i \right ] .
\end{array}
\end{equation}
which is exactly the same as the one presented by Wang et al. in (\cite{wang}, Example 4.1, equations (61)-(62)).

\subsection{Variational integrators}

Classical results on order one embedding ensures that critical points of the discrete functional $\mathscr{L}_{H_{W-Z,\epsilon},h,1}$ are given by (see \cite{rouba1}):
\begin{align}
\label{order1}
\Delta_{-} [p] &=-\di\frac{\partial H}{\partial q} (p,q) ,\\
\Delta_{+} [q]& =\di\frac{\partial H}{\partial p} (p, q) ,
\end{align}
for all $t\in \T^{\pm}$.\\

Here again, we recover the algebraic form of classical Hamiltonian systems. Indeed, the structure of these equations can be written as follows in the classical case:

\begin{equation}
    \left .
    \begin{array}{lll}
    -\left ( \di\frac{d}{dt} \right )_{*} (p) & = & -\partial_q H (p,q),\\
    \di\frac{d}{dt} & = & \partial_p H (p,q) ,
    \end{array}
    \right .
\end{equation}
where $\left ( \di\frac{d}{dt} \right )_{*} = -\di\frac{d}{dt}$ stands for the adjoint differential operator associated to $\di\frac{d}{dt}$ for the usual scalar product $\langle f,g\rangle = \di\int_{\R} fg dt$. \\

In the discrete case, the adjoint of $\Delta_{+}$ with respect to the discrete scalar product $(f,g)=\di\int_a^b f_{\circ} \Delta_{+} (g) \Delta t$ is given by  $-\Delta_{-}$ so that if $O$ stands for the differential operator 
\begin{equation}
    O \left ( \di\frac{d}{dt} \right ) =
    \di\left ( 
    \left . 
    \begin{array}{c}
    -\left ( \di\frac{d}{dt} \right )_{*} \\
    \di\frac{d}{dt}
    \end{array}
    \right .
    \right ) 
\end{equation}
acting on vector of functions $(p,q)$, the classical Hamiltonian system can be written as 
\begin{equation}
\label{hamnew}
    O \left ( \di\frac{d}{dt} \right )\left [
    \begin{array}{l}
    p\\
    q
    \end{array}
    \right ]
    = J \Delta H (p,q) ,    
\end{equation}
where $\Delta H$ is the gradient vector $\Delta H = \left ( \begin{array}{l}
\di\frac{\partial H}{\partial p}\\
\di\frac{\partial H}{\partial q}
\end{array}
\right )$ and $J$ is the simplectic matrix $J=\left ( 
\begin{array}{cc}
0 & -I_d \\
I_d & 0 
\end{array}
\right )$ where $I_d$ is the $d\times d$ identity matrix. \\

We see directly that the discrete equation \eqref{order1} can be written as 
\begin{equation}
    O (\Delta_+ )\left [
    \begin{array}{l}
    p\\
    q
    \end{array}
    \right ]
    = J \Delta H (p,q) ,   
\end{equation}
which possesses the same algebraic form as \eqref{hamnew}.

We have proved in \cite{rouba2} that critical points of the Wong-Zakai discrete functional \eqref{wong-zakai-discrete} corresponds to solutions of the discrete Hamiltonian system
\begin{align}
\Delta_{1/2,-} [p_{\circ}] &=\left [ -\di\frac{\partial H}{\partial q} (p_{\circ},q_{\circ}] \right ]_{1/2,-} ,\\
\Delta_{\circ,+} [q_{\circ}]& =\di\frac{\partial H}{\partial p} (p_{\circ}, q_{\circ} ) ,
\end{align}
for all $t\in \T^{\pm}_{1/2}$.\\

\section{Numerical examples}
\label{example}
\subsection{The stochastic Kubo oscillator}

The problem of convergence of Wong-Zakai variational integrators to the solutions of the stochastic Hamiltonian system can be illustrated in the Kubo case. As we have said, the fact that Stratonovich stochastic integrals are defined using mid-point quadrature formula implies that using only order one quadrature for the Wong-Zalai approximation will lead to wrong simulations. As an example, taking as initial conditions $p_0=1$, $q_0=0$ and $h=\epsilon=0.02$, for $a=1.5$ and $b=1$, we obtain

\begin{figure}[htb!]
    \includegraphics[width=0.4\textwidth]{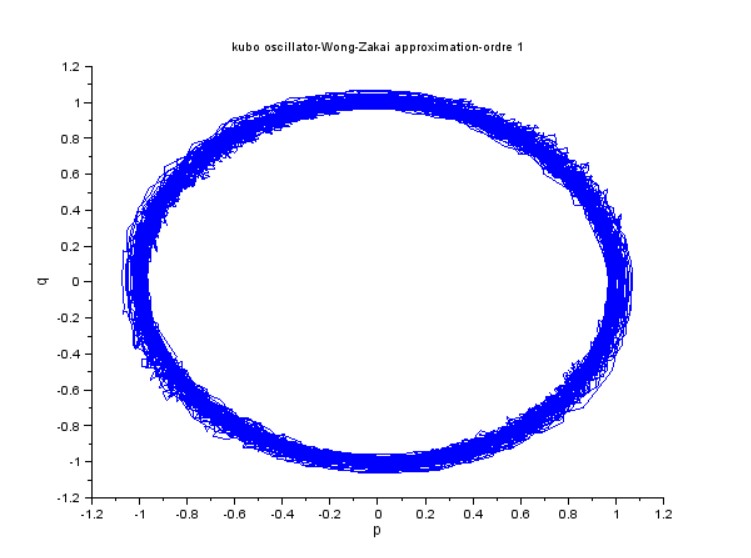}
    \includegraphics[width=0.4\textwidth]{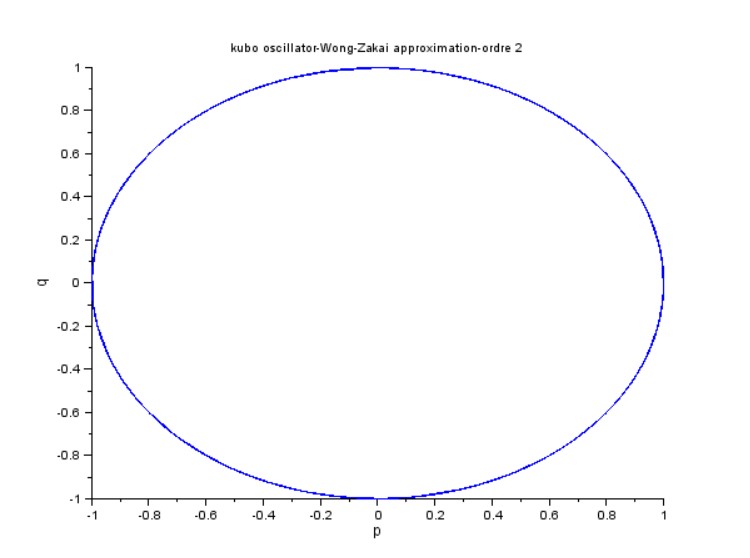}
    \caption{Order 1 (left) and order 2 (right) Wong-Zakai variational integrators}
    \label{fig1}
\end{figure}

As expected the case of order $1$ does not give a satisfying result instead of the mid-point one. \\

We can also look for the behaviour of the Wong-Zakai energy and of the unperturbed energy over Wong-Zakai simulations with the two variational integrators. 

\begin{figure}[ht!]
    \includegraphics[width=0.4\textwidth]{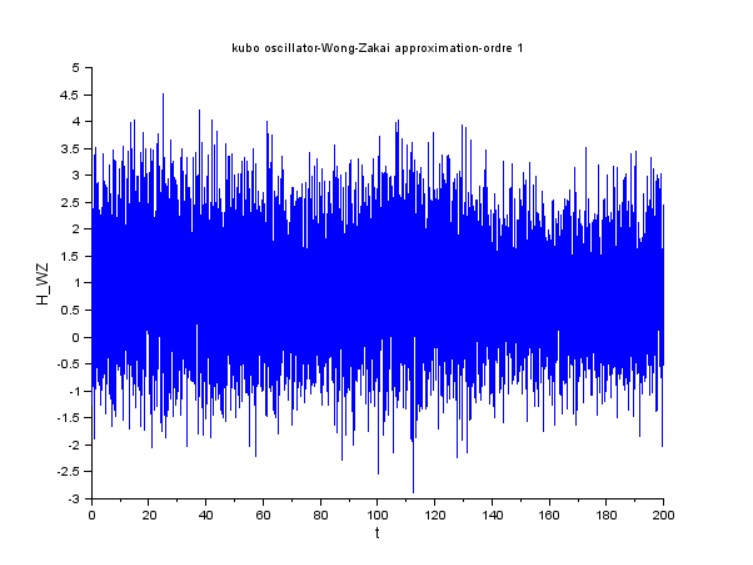}
    \includegraphics[width=0.4\textwidth]{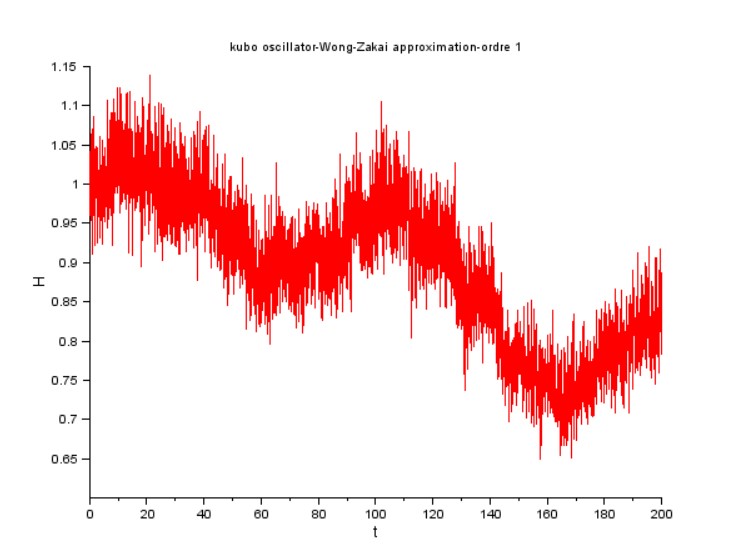}
    \caption{Wong-Zakai energy $H_{WZ, \epsilon}$ (left)  and unperturbed energy $H$ (right) with an order 1 Wong-Zakai variational integrators}
    \label{fig2}
\end{figure}

In figure \ref{fig2}, we see that the unperturbed energy is not preserved indicating that the order one Wong-Zakai integrator is not appropriate. 

\begin{figure}[ht!]
    \includegraphics[width=0.4\textwidth]{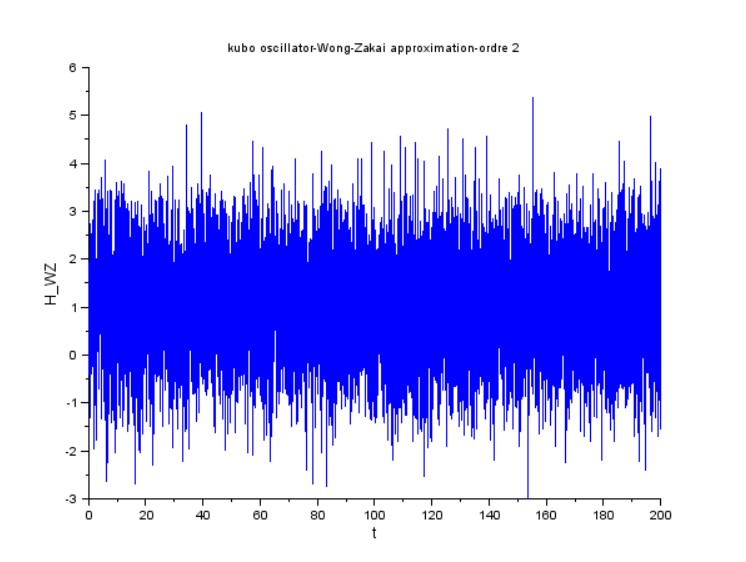}
    \includegraphics[width=0.4\textwidth]{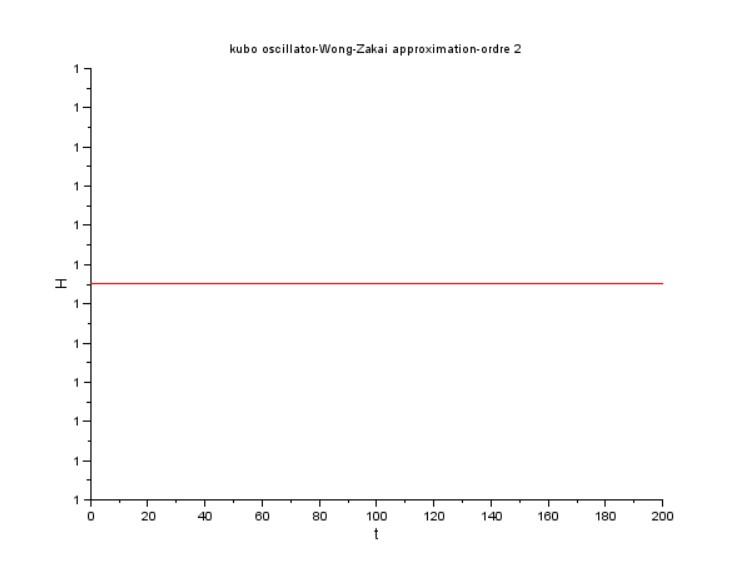}
    \caption{Wong-Zakai energy $H_{WZ, \epsilon}$ (left) and unperturbed energy $H$ (right) with an order 2 Wong-Zakai variational integrators}
    \label{fig3}
\end{figure}

Figure \ref{fig3} shows that the mid-point variational integrator very well behaves with respect to the preservation of the unperturbed energy. We see also that the behavior of the Wong-Zakai energy does not give many insights on the dynamics of the stochastic Hamiltonian system.

\subsection{Hamiltonian systems with two additive noises}
Let us consider the stochastic Hamiltonian system with two additive noises defined by
\begin{align}
\label{ex2-1}
    dp & = -q dt + \gamma \circ dW_2(t),\\
    \label{ex2-2}
    dq & = pdt + \sigma \circ dW_1(t),
\end{align}
with 
\begin{align}
    H(p, q)=\frac{1}{2} (p^2 +q^2),&& H_1(p,q)= \sigma p ,&& H_2(p,q) = -\gamma q.
\end{align}
Note that for systems with additive noise, It\^o and stratonovich are equivalent due to the fact that the Wang-Zakai correction vanishes when the drift coefficients are constant.\\

The  Wong-Zakai variational integrators of order one for the system \eqref{ex2-1}-\eqref{ex2-2} are given by 
\begin{align}
\label{Ex2-VI1}
    p_{i} &= p_{i-1} - h q_{i} + h \gamma \dot{W}_{eps,2,i}, \\
\label{Ex2-VI2}
    q_{i+1}& = q_{i} - h p_{i} + h \sigma \dot{W}_{eps,1,i},
\end{align}
we consider the following reference solution  for the system \eqref{ex2-1}-\eqref{ex2-2}  defined by (see \cite{Milstein})
\begin{equation}
    \begin{array}{ll}
       \tilde{X}_{i+1} & = \tilde{H} \tilde{X}_i + \tilde{v}_k  \\
         &= \begin{bmatrix}
             cos\ h & sin \ h\\
             -sin \ h & cos \ h
         \end{bmatrix} 
         \tilde{X}_i + \di\frac{1}{h}
         \begin{bmatrix}
             \sigma sin\ h\, \Delta_i W_1 + 2 \gamma sin^2 \di\frac{h}{2} \, \Delta_i W_2\\
             -2 \sigma sin^2 \di\frac{h}{2}\,  \Delta_i W_1 + \gamma sin\ h \,\Delta_i W_2 
         \end{bmatrix}
         ,
    \end{array}
\end{equation}
where $ \tilde{X}_{i} = (\tilde{p}_i, \tilde{q}_i)$.\\

\begin{figure}[ht!]
    \includegraphics[width=0.7\textwidth]{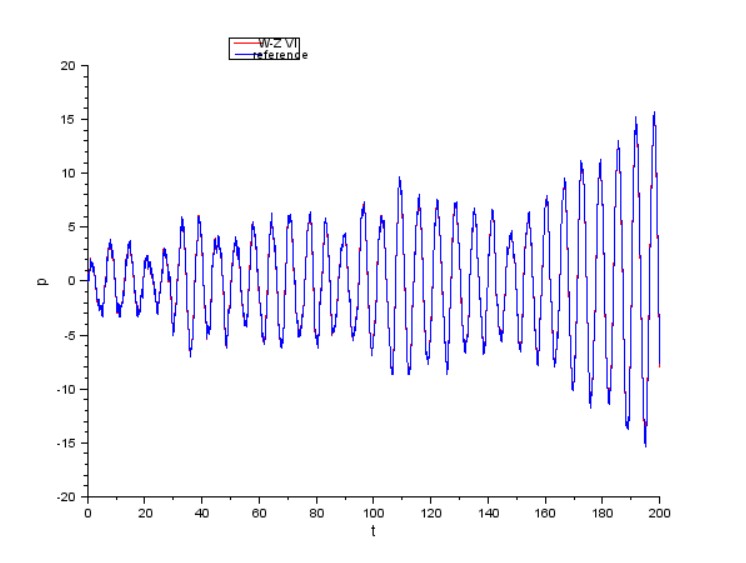}
    \caption{Order $1$ Wong-Zakai variational integrators (red) and the reference solution (blue).}
    \label{fig4}
\end{figure}

In figure \ref{fig4}, we present the simulations for the reference solution and the numerical integrators obtained using the  Wong-Zakai variational integrators \eqref{Ex2-VI1}-\eqref{Ex2-VI2}, with initial conditions $p_0=0$ and $q_0=0$ for $\sigma=0$, $\gamma=1$, $h=0.02$, $\epsilon=0.02$ and $t \le 200$.

As we see in figure \ref{fig4}, the two paths coincide visually. We can deduce that in the case of systems with additive noise, an order $1$ quadrature formula for both deterministic and random part provides a good approximation for the exact solution  seeing the fact that the stratonovich setting is not important anymore meaning that a quadrature of order $2$ for the random part  is not a necessary conditions to have convergence.

\section{Conclusion and perspectives}
In this article, we have given a rigorous meaning to computations made in Wang et al. \cite{wang} using Wong-Zakai approximation of diffusion and discrete embedding point of view. As expected, the convergence of the resulting numerical scheme toward the solutions of the stochastic Hamiltonian system necessitates quadrature formula of order two of the random part at least when the noise is not additive. In order to be complete, one need to compare these variational integrators with the one obtained directly from the stochastic variational formulation of stochastic Hamiltonian systems. Constructions of such integrator were already made by Bou-Rabee et al. in \cite{bou-rabee} for a particular class of stochastic Hamiltonian systems using the Marsden's strategy. Here again, a stochastic discrete embedding can be used to reinterpret these constructions. This will be done in a forthcoming paper.

\end{document}